# One Categorization of Microtonal Scales

Luka Milinković[a], Branko Malešević[b], Dragana Pavlović –Šumarac[a], Bojan Banjac[b,c], Miomir Mijić[a]

[a] *Department of Telecommunication, School of Electrical Engineering, University of Belgrade, Belgrade, Serbia*

[b] *Department of Applied Mathematics, School of Electrical Engineering, University of Belgrade, Belgrade, Serbia*

[c] *Department of Fundamentals Sciences, Faculty of Technical Sciences, University of Novi Sad, Novi Sad, Serbia*

Corresponding author is Branko Malešević (malesevic@etf.bg.ac.rs)

Luka Milinković (lukaui@yahoo.com), Branko Malešević (malesevic@etf.bg.ac.rs), Dragana Pavlović –Šumarac (dsumarac@etf.bg.ac.rs), Bojan Banjac (bojan.banjac@uns.ac.rs), Miomir Mijić (emijic@etf.bg.ac.rs)

# One Categorization of Microtonal Scales

This study considers rational approximations of musical constant $\beta = log_2\left(\frac{3}{2}\right)$, which defines perfect fifth. This constant has been the subject of the numerous studies, and this paper determines quality of rational approximations in regards to absolute error. We analysed convergents and secondary convergents (some of these are the best Huygens approximations). Especially, we determined quality of the secondary convergents which are not the best Huygens approximations – in this paper we called them non-convergents approximations. Some of the microtonal scales have been positioned and determined by using non-convergents approximation of music constant which defines perfect fifth.

## Introduction

Continued fractions belong to the field of mathematics that has influence on theory of music. Much about this topic has been written in 20th century. At that time, convergents and secondary convergents were considered as the part of the continued fraction theory by mathematician A. Ya. Khinchin (Khinchin 1964). We have decided to perform further research this field, and to check whether some fractions with similar optimal properties to convergents and secondary convergents of a continued fraction exists, and to consider their applications.

One of the many important applications of the convergents and secondary convergents is in music theory. The most important intervals in music are octave, perfect fifth, and major third, which naturally fit to our hearing. They are based on simple relations 2:1, 3:2, and 5:4, respectively. Perfect fifth and major third have been analysed into more details in papers (Benson D. 2006; Kent J. T. 1986; Krushchev S. 2008.), while in the paper (Liern V. 2015) has been noted through the example 5.9, specific array of fractions which approximates music constant, which in turn defines perfect fifth

$$\beta = log_2\left(\frac{3}{2}\right) = 0{,}5849625007211\ldots \tag{1}$$

Rational approximations of the music constant $\beta$, which are discussed in the paper (Liern V. 2015) in the example 5.9, are continued fractions. The goal of this study is to show that there are other approximations of the previously mentioned constant. In this paper, we will do the following:

(1) Determine rational approximations of the first kind (Khinchin 1964), which are not continued fractions of the constant β

(2) Consider possibilities of categorization of remaining rational approximations which are not the best rational approximations of the first kind of the constant β.

Considerations of the problem 1) and 2) will solve the general problem of categorization of rational approximations for arbitrary real number. Also, the application of categorization will be considered, with the goal to determine microtonal scales. After the analysis of the problems 1) and 2) in the study, pseudo code which determines quality of rational approximation is given.

## Continued Fractions

Continued fraction is expression given in the following form:

$$\alpha = a_0 + \cfrac{1}{a_1 + \cfrac{1}{a_2 + \cfrac{1}{a_3 + \ddots}}} \tag{2}$$

where $\alpha \in \mathbb{N}$, $a_0 \in \mathbb{N}_0 = \mathbb{N} \cup \{0\}$, $a_i \in \mathbb{N}(i \geq 1)$. Throughout the paper we will consider only cases when α is positive number, while whole analysis can be done analogously in case when α is negative number. Continued fractions are written in shorter form $[a_0; a_1, a_2, \ldots, a_n]$. For convergents

$$\frac{p_n}{q_n} = [a_0; a_1, a_2, \ldots, a_{n-1,} a_n] \tag{3}$$

we introduce *secondary convergents*

$$\frac{p'_n}{q'_n} = [a_0; a_1, a_2, \ldots, a_{n-1}, a'_n] \tag{4}$$

where $a'_n \in \mathbb{N} \wedge 0 < a'_n < a_n$.

**Properties of Continued Fractions**

Basic properties of continued fractions are described through the following.

**Property 1.** $\mathrm{GCD}(p_n, \mathrm{q_n}) = 1;\ (n \in \mathbb{N}_0) \wedge q_n \nearrow^{\infty}\ (n \longrightarrow \infty)$

**Property 2.** $\lim_{n\longrightarrow\infty} \left|\frac{p_{n+1}}{q_{n+1}} - \frac{p_n}{q_n}\right| = \lim_{n\longrightarrow\infty} \frac{1}{q_n \cdot q_{n+1}} = 0$

**Property 3.** $\frac{p_0}{q_0} < \frac{p_2}{q_2} < \cdots \leq \alpha \leq \cdots < \frac{p_3}{q_3} < \frac{p_1}{q_1};\ (\alpha \in \mathbb{R})$

**Property 4.** $\left|\alpha - \frac{p_n}{q_n}\right| < \frac{1}{q_{n-1}^2};\ (\alpha \in \mathbb{R} \wedge n \in \mathbb{N})$

For the fraction $\frac{p}{q}$ value of *absolute error* Δ in approximation of the constant α is given by

$$\Delta = \left|\alpha - \frac{p}{q}\right|. \tag{5}$$

**Definition 1.** For $\alpha \in \mathbb{R}$, the fraction $\frac{p}{q}$ is *approximation of the first kind* (the best Huygens approximation) if for every fractions $\frac{r}{s}$ stands $(0 < s \leq q)$:

$$\left|\alpha - \frac{p}{q}\right| < \left|\alpha - \frac{r}{s}\right|. \tag{6}$$

In the book (Krushchev S. 2008) approximations of the first kind are called *the best Huygens approximations*.

**Definition 2.** For $\alpha \in \mathbb{R}$, the fraction $\frac{p}{q}$ is *approximation of the second kind* if for every fraction $\frac{r}{s}$ stands $(0 < s \leq q)$:

$$|q\alpha - p| < |s\alpha - r|. \quad (7)$$

**Property 5** (Khinchin 1964). Every best approximation of the second kind is convergent. Every best approximation of the first kind can be convergent or secondary convergent.

**Remark.** Sequence of approximations of the second kind consists only of convergents. We called them *continued approximations* (*continued fractions*) and labeled them as (c). Besides them, array of approximations of the first kind can contain some secondary convergents, and those fractions are called *semi-continued approximations* (*semi-continued fractions*) and are labeled with (sc). A. Ya. Khinchin called semi-continued approximations *intermediate fractions* [1].

Secondary convergents $\frac{p_n'}{q_n'}$ are considered *semi-continued approximation* if following is true

$$\left|\alpha - \frac{p_n}{q_n}\right| < \left|\alpha - \frac{p_n'}{q_n'}\right| < \left|\alpha - \frac{p_{n-1}}{q_{n-1}}\right|. \quad (8)$$

Determination of convergents and secondary convergents can be done in one of the ways described in the papers (Cortzen A. 2011; Malešević B. 1988; Malešević B, Milinković L. 2014.).

**Property 6.** If there exist semi-continued approximation$[a_0; a_1, \ldots, a_{n-1}, a_n']$, with lowest value $a_n'$, then fractions $[a_0; a_1, \ldots, a_{n-1,} a_n' + j]$, for $0 < j < a_n - a_n'$, are semi-continued approximations too.

Properties 1-5 are described based on (Khinchin 1964), while Property 6 is based on (Malešević B. 1988.)

**Example 1.** We will use number $\pi$ to show that between two continued fractions can be found fractions which can be better or worse approximation of number $\pi$ when compared with existing continued fractions. Continued fraction of the number $\pi$ is

$$\pi = 3.14159265359 \ldots = [3; 7,15,1,292,1,1,1,2,1,1,3 \ldots]. \qquad (9)$$

Number $\pi$ and his initial rational approximations are shown on the Figure 1. First continued approximation is $\frac{3}{1}$ (A), and second continued approximation is $\frac{22}{7}$ (G). First semi-continued approximation is $\frac{13}{4}$ (D), second semi-continued approximation is $\frac{16}{5}$ (E), and third semi-continued approximation is $\frac{19}{6}$ (F).

There are another two rational approximations $\frac{6}{2}$ and $\frac{9}{3}$ (B and C) between two consecutive continued approximations. For denominator 2 numerator is determined as $6 = \lfloor 2\pi \rfloor$, while for denominator 3 numerator is determined as 9= $\lfloor 3\pi \rfloor$. These two fractions are not better approximations of number $\pi$ then first fraction, $\frac{3}{1}$. They are duplicates of $\frac{3}{1}$. These two fractions in section *Defining Categories* are called multiples.

Figure 1 is showing first continued fraction (A) that represents the worst approximation (with the biggest absolute error) of number $\pi$. Every next semi-continued or continued approximation (from D to G) has diminishing absolute error while tends to be closer towards real value of $\pi$.

Please note, on Figure 1 the values of fraction are shown on ordinate and the values of denominators of fractions are shown on abscissa.

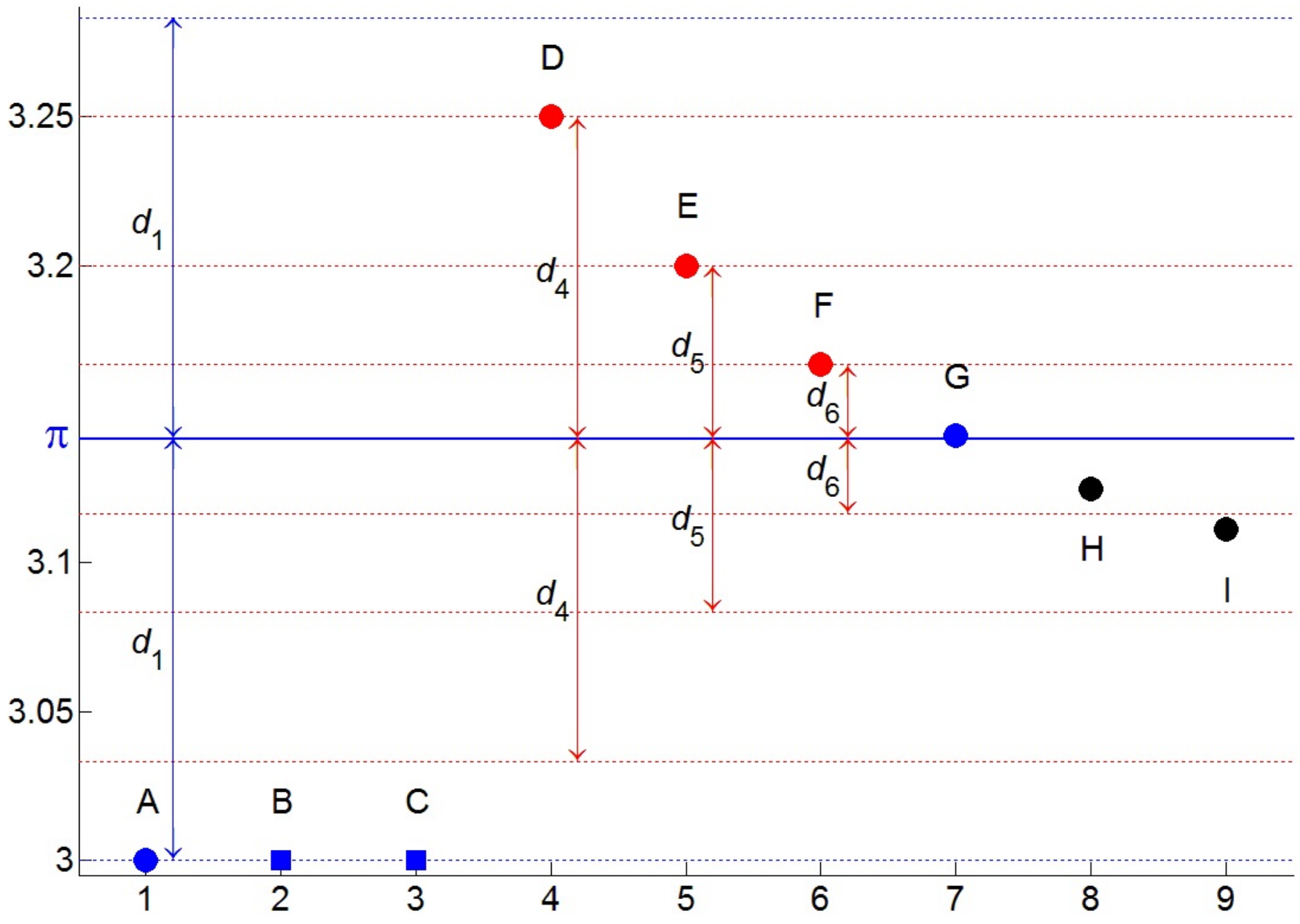

Figure 1. Rational approximation of the number $\pi$

In this paper, fractions that are neither continued nor semi-continued are called *non-continued approximations* (*non-continued fractions*). Based on previous consideration, non-continued fractions are second convergents which are not the best Huygens approximations. For this example, non-continued fractions are $\frac{25}{8}$ and $\frac{28}{9}$ (H and I), Figure 1.

■

Convergents and secondary convergents are being considered and applied in many different fields, among others in following:

- Music (Barbour J. M. 1948; Benson D. 2006; Carey N, Clampitt D. 1989; Douthett J, Krantz R. 2007; Hall R. W, Josi K. 2001; Haluska J. 2003; Kent J. T. 1986; Krantz R. J. 1998; Krantz R. J, Douthett J. 2001; Krushchev S. 2008;

Liern V. 2015; Noll T. 2006; Savage M 2014; Schechter M. 1980; Smoyer L. 2005)

- Astronomy for approximation of various constants related to universe, such as tropical year, and synodic moon (Malešević B, Milinković L. 2014),
- Mathematics for approximation of math constants such as Euler's number e, then π and others (Malešević B. 1988; Cortzen A. 2011.)
- Mechanic at Antikythera Mechanism (Malešević B, Milinković L. 2014; Popkonstantinović P, Obradović R, Jeli Z, Mišić S. 2014;),
- Cryptography at RSA algorithms and in generation of random numbers (Malešević B, Milinković L. 2012) etc.

**Non-continued Fractions and Their Properties**

In this paper we consider the best rational approximations of the real number α in the following form

$$\frac{k_i}{m_i} \tag{10}$$

where $k_i, m_i \in \mathbb{N}$ and denominator are defined with $m_i = i$ for $i = 1,2,\ldots,\mathbb{N}$.

Now we form an array

$$\boldsymbol{A} = \left(\frac{k_i}{m_i} :\ i = 1,2,\ldots\right) \tag{11}$$

where

$$m_i = i \quad \wedge \quad k_i = \mathbf{round}(m_i) = \begin{cases} \lfloor \alpha \cdot m_i + 0{,}5 \rfloor, & \lfloor \alpha \cdot m_i \rfloor \neq 2t \\ \lfloor \alpha \cdot m_i \rfloor, & \lfloor \alpha \cdot m_i \rfloor = 2t \end{cases} \tag{12}$$

for some $t \in \mathbb{N}$. On this way we choose numerator $k_i$ for given denominator $m_i$ with minimal absolute error $\Delta_i$.

**Property 7.** Every fraction from array $\boldsymbol{A}$ that is better approximates value $\alpha$ than its first previous determined best approximation is also continued or semi-continued fraction [3].

The music constant $\beta$ (Eq. 1) has the following continued fraction form

$$\beta = [0; 1,1,2,2,3,1,5,2,2,3,2,2,1,1,5,5 \dots], \tag{13}$$

and array $\boldsymbol{A}_\beta$ is given with

$$\boldsymbol{A}_\beta = \Big(\frac{1}{1}(\text{c}), \frac{1}{2}(\text{c}), \frac{2}{3}(\text{sc}), \frac{2}{4}, \frac{3}{5}(\text{c}), \frac{4}{6}, \frac{4}{7}(\text{sc}), \frac{5}{8}, \frac{5}{9}, \frac{6}{10}, \frac{7}{11}, \frac{7}{12}(\text{c}),$$
$$\frac{8}{13}, \frac{8}{14}, \frac{9}{15}, \frac{9}{16}, \frac{10}{17}, \frac{11}{18}, \frac{11}{19}, \frac{12}{20}, \frac{12}{21}, \frac{13}{22}, \frac{13}{23}, \frac{14}{14}, \frac{15}{25}, \frac{15}{26}, \frac{16}{27}, \frac{16}{28}, \tag{14}$$
$$\frac{17}{29}(\text{sc}), \frac{18}{30}, \frac{18}{31}, \frac{19}{32}, \frac{19}{33}, \frac{20}{34}, \frac{20}{35}, \frac{21}{36}, \frac{22}{37}, \frac{22}{38}, \frac{23}{39}, \frac{23}{40}, \frac{24}{41}(\text{c}) \dots\Big)$$

Properties which will further be considered are new.

**Property 8.** Let us consider the array from $\frac{k_i}{m_i}$ to $\frac{k_j}{m_j}$, where $i + 1 < j$ (minimum 3 elements), such that only $\frac{k_i}{m_i}$ and $\frac{k_j}{m_j}$ are continued or semi-continued fractions. Then the following inequalities are true

$$\Delta_j < \Delta_i < \Delta_l \tag{15}$$

where $i \leq l \leq j$.

Property 8 is described with following example.

**Example 2.** In array $\boldsymbol{A}_\beta$, for the constant $\beta$, between semi-continued fraction $\frac{4}{7}(\text{sc})$ and continued fraction $\frac{7}{12}(\text{c})$ exist fractions $\frac{5}{8}$, $\frac{5}{9}$, $\frac{6}{10}$ and $\frac{7}{11}$. Absolute error values are

$$\Delta_7(\text{sc}) = 0{,}013533 \dots \tag{16}$$

$$\Delta_8 = 0{,}040037 \dots \tag{17}$$

$$\varDelta_9 = 0{,}029406\ldots \quad (18)$$

$$\varDelta_{10} = 0{,}015037\ldots \quad (19)$$

$$\varDelta_{11} = 0{,}039507\ldots \quad (20)$$

$$\varDelta_{12}(\mathrm{c}) = 0{,}001629\ldots \quad (21)$$

We can notice that the biggest absolute errors are $\varDelta_8, \ldots, \varDelta_{11}$, according the Property 8.

■

**Property 9.** For real number α exist non-continued approximations $\frac{k_i}{m_i}$, which are better approximations than continued or semi-continued approximations $\frac{k_j}{m_j}$, where $m_j < m_i$. Thereby, non-continued approximation $\frac{k_i}{m_i}$ is never better than first previous continued or semi-continued approximation $\frac{k_l}{m_l}$, whereby $m_j < m_l < m_i$. Than it is

$$\varDelta_j > \varDelta_i > \varDelta_l. \quad (22)$$

Property 9 is illustrated on example for the constant $\beta$.

**Example 3.** How it is shown in previous example, non-continued fraction $\frac{5}{8}$ is not a better approximation of the constant $\beta$ than its first previous semi-continued fraction $\frac{4}{7}(\mathrm{sc})$. However, $\frac{5}{8}$ is a better approximation of the constant $\beta$ than semi-continued fraction $\frac{2}{3}(\mathrm{sc})$, from array $\boldsymbol{A}_\beta$, with absolute error

$$\varDelta_3(\mathrm{sc}) = 0{,}081704\ldots \quad (23)$$

From this for

$$m_3 = 3 < m_7 = 7 < m_8 = 8 \quad (24)$$

stands that

$$\varDelta_3 = 0{,}081704 > \varDelta_8 = 0{,}040037 > \varDelta_7 = 0{,}01533. \quad (25)$$

The same conclusion is true for the rest non-continued fractions from Example 2 $\left(\frac{5}{9} \text{ and } \frac{7}{11}\right)$.

■

**Property 10.** Let us consider two fractions $\frac{k_i}{m_i}$ and $\frac{k_j}{m_j}$. Let $\frac{k_i}{m_i}$ be non-continued and $\frac{k_j}{m_j}$ be continued or semi-continued fraction. Let for them hold true

$$\varDelta_j > \varDelta_i. \quad (26)$$

Than $m_j < m_i$ is true. For all fractions $\frac{k_l}{m_l}$, where $0 < m_l < m_j$, we have

$$\varDelta_l > \varDelta_i. \quad (27)$$

Property 10 is illustrated on example for the constant $\beta$.

**Example 4.** For constant $\beta$, which is considered in the Examples 3, can be also proven that continued fractions $\frac{1}{1}$ (c) and $\frac{1}{2}$ (c), which are positioned in array $\boldsymbol{A}_\beta$ before semi-continued fraction $\frac{2}{3}$ (sc), are worse approximations of the constant $\beta$ than non-continued fraction $\frac{5}{8}$. Since

$$\varDelta_3 = 0{,}081704 > \varDelta_8 = 0{,}040037, \quad (28)$$

than

$$\varDelta_2 = 0{,}084962 > \varDelta_8 = 0{,}040037 \quad (29)$$

and

$$\Delta_1 = 0{,}415037 > \Delta_8 = 0{,}040037. \quad (30)$$

■

**Defining Categories**

Let us observe a sub-array ***B***, of the array ***A***, which is composed only of continued and semi-continued approximations (without non-continued). In this sub-array, it is clear that every next fraction is better approximation of value α than previous approximation [4]. In general case, when we have continued, semi-continued and non-continued fractions in the same array, previous conclusion is not true. Therefore, it is necessary to determine category for each fraction in regard to the value of absolute error.

For rational values α, array ***A*** is finite and can contain large number of members. On the other hand, for irrational values of α, array ***A*** always has infinite number of members. In practical applications, assigning categories is done over the finite sub-array of the array ***A***. In this article, we will mark ***B*** as finite sub-array of the array ***A*** and define it as

$$\boldsymbol{B}\left(\frac{k_s}{m_s}, \frac{k_e}{m_e}\right) = \left(\frac{k_s}{m_s}, \frac{k_{s+1}}{m_{s+1}}, \frac{k_{s+2}}{m_{s+2}}, \dots, \frac{k_{e-2}}{m_{e-2}}, \frac{k_{e-1}}{m_{e-1}}, \frac{k_e}{m_e}\right), \quad (31)$$

where $m_s < m_e$. First element of array ***B*** has index *s*, and the last one index *e*. In accordance with the above, number of elements or array ***B*** is determined with

$$M = e - s + 1. \quad (32)$$

In the array ***B***, we define members of 1st category as fractions that have the best approximation of value $\alpha$. Among the members of the array ***B***, which have not fallen in to the 1st category, those that have the best approximations of value $\alpha$ are defined as fractions of 2nd category. Among the rest of the fractions in array $\boldsymbol{B}$ we define category by induction.

All properties specified in this paper (Properties 1-10), which can be applied on whole array $\boldsymbol{A}$, can also be applied to any of its sub-arrays.

Let us notice that in array $\boldsymbol{B}\left(\frac{k_s}{m_s}, \frac{k_e}{m_e}\right)$ we can select following sub-arrays:

- $\boldsymbol{B^c}\left(\frac{k_s}{m_s}, \frac{k_e}{m_e}\right)$ – array of continued fractions and their multiples in $\boldsymbol{B}\left(\frac{k_s}{m_s}, \frac{k_e}{m_e}\right)$,
- $\boldsymbol{B^{sc}}\left(\frac{k_s}{m_s}, \frac{k_e}{m_e}\right)$ – array of semi-continued fractions and their multiples in $\boldsymbol{B}\left(\frac{k_s}{m_s}, \frac{k_e}{m_e}\right)$,
- $\boldsymbol{B^v}\left(\frac{k_s}{m_s}, \frac{k_e}{m_e}\right)$ – array of non-continued fractions and their multiples in $\boldsymbol{B}\left(\frac{k_s}{m_s}, \frac{k_e}{m_e}\right)$.

Consequently, for the following sub-arrays we can say

$$\left\{\boldsymbol{B^c}\left(\frac{k_s}{m_s}, \frac{k_e}{m_e}\right)\right\} \cap \left\{\boldsymbol{B^{sc}}\left(\frac{k_s}{m_s}, \frac{k_e}{m_e}\right)\right\} = \emptyset, \qquad (33)$$

$$\left\{\boldsymbol{B^c}\left(\frac{k_s}{m_s}, \frac{k_e}{m_e}\right)\right\} \cap \left\{\boldsymbol{B^v}\left(\frac{k_s}{m_s}, \frac{k_e}{m_e}\right)\right\} = \emptyset, \qquad (34)$$

$$\left\{\boldsymbol{B^{sc}}\left(\frac{k_s}{m_s}, \frac{k_e}{m_e}\right)\right\} \cap \left\{\boldsymbol{B^v}\left(\frac{k_s}{m_s}, \frac{k_e}{m_e}\right)\right\} = \emptyset, \qquad (35)$$

and

$$\left\{\boldsymbol{B}\left(\frac{k_s}{m_s}, \frac{k_e}{m_e}\right)\right\} = \left\{\boldsymbol{B^c}\left(\frac{k_s}{m_s}, \frac{k_e}{m_e}\right)\right\} \cup \left\{\boldsymbol{B^\mu}\left(\frac{k_s}{m_s}, \frac{k_e}{m_e}\right)\right\} \cup \left\{\boldsymbol{B^v}\left(\frac{k_s}{m_s}, \frac{k_e}{m_e}\right)\right\}. \qquad (36)$$

Let us define two categorizations for array $\boldsymbol{B}$, (Eq. 31):

(1) <u>Total categorization ($\tau$)</u>

For each fraction in array $\boldsymbol{B}\left(\frac{k_s}{m_s},\frac{k_e}{m_e}\right)$ we will calculate category based on absolute error starting from 1st *total category*. Element from this array which has the smallest absolute error of value α belongs to 1st (total) category. The result of this sorting of array $\boldsymbol{B}$ is array $\boldsymbol{\mathcal{B}}^{\boldsymbol{\tau}}\left(\frac{k_s}{m_s},\frac{k_e}{m_e}\right)$.

(2) Non-continued categorization (υ)

For each fraction in array $\boldsymbol{B}^{\boldsymbol{\upsilon}}\left(\frac{k_s}{m_s},\frac{k_e}{m_e}\right)$ we will calculate category based on absolute error start from 1st *non-continued category*. Element from this array which has the smallest absolute error of value α belongs to 1st (non-continued) category. The result of this sorting of array $\boldsymbol{B}^{\boldsymbol{\upsilon}}$ is array $\boldsymbol{\mathcal{B}}^{\boldsymbol{\upsilon}}\left(\frac{k_s}{m_s},\frac{k_e}{m_e}\right)$.

Let us emphasize that with $\boldsymbol{\mathcal{B}}$ is marked sorted array $\boldsymbol{B}$. Therefore, for both, total and non-continued categorization, multiples will be of the exactly same categories as their basic fractions, but specifically marked as it will be shown.

In order to explain the preceding definitions on example of the constant $\beta$ we will form array $\boldsymbol{B}_\beta$ based on 30 elements of array $\boldsymbol{A}_\beta$ (Eq. 14). Range was formed based on common musical scale $\left(\frac{7}{12}\right)$ and two consecutive continued fractions.

$$\boldsymbol{B}_\beta\left(\frac{7}{12},\frac{24}{41}\right)=\left(\frac{7}{12},\frac{8}{13},\frac{8}{14},\frac{9}{15},\frac{9}{16},\frac{10}{17},\frac{11}{18},\frac{11}{19},\frac{12}{20},\frac{12}{21},\frac{13}{22},\frac{13}{23},\frac{14}{24},\frac{15}{25},\frac{15}{26},\right. \tag{37}$$
$$\left.\frac{16}{27},\frac{16}{28},\frac{17}{29},\frac{18}{30},\frac{18}{31},\frac{19}{32},\frac{19}{33},\frac{20}{34},\frac{20}{35},\frac{21}{36},\frac{22}{37},\frac{22}{38},\frac{23}{39},\frac{23}{40},\frac{24}{41}\right).$$

In order to form array of non-continued fractions $\boldsymbol{B}_\beta^{\boldsymbol{\upsilon}}$, it is necessary to omit continued and semi-continued fractions of the constant $\beta$, and multiples of continued and semi-continued fractions from array$\boldsymbol{B}_\beta$. Rejected elements are $\frac{7}{12},\frac{17}{29},\frac{24}{41}$, as well as multiples $\frac{8}{14},\frac{9}{15},\frac{12}{20},\frac{12}{21},\frac{14}{24},\frac{15}{25},\frac{16}{28},\frac{18}{30},\frac{20}{35},\frac{21}{36}$. Based on the mentioned properties we find continued and semi-continued fractions and their multiples in the array $\boldsymbol{B}_\beta$. The resulting array $\boldsymbol{B}_\beta^{\boldsymbol{v}}$ is

$$\boldsymbol{B}_\beta^{\boldsymbol{v}}\left(\frac{7}{12},\frac{24}{41}\right)=\left(\frac{8}{13},\frac{9}{16},\frac{10}{17},\frac{11}{18},\frac{11}{19},\frac{13}{22},\frac{13}{23},\frac{15}{26},\frac{16}{27},\frac{18}{31},\frac{19}{32},\frac{19}{33},\frac{20}{34},\frac{22}{37},\frac{22}{38},\frac{23}{39},\frac{23}{40}\right). \quad (38)$$

The best approximation in the array $\boldsymbol{B}_\beta$ is the last continued or semi-continued fraction (Kent J. T. 1986; Khinchin 1964; Krushchev S. 2008.). In the observed array, this is a continued fraction $\frac{24}{41}$. If the array had been shorter, for example the last element could have been $\frac{23}{40}$, than the best approximation would be given semi-continued fraction $\frac{17}{29}$, which is, in this observed array $\boldsymbol{B}_\beta$, second in total categorization.

Table 1 shows all elements of the array $\boldsymbol{B}_\beta\left(\frac{7}{12},\frac{24}{41}\right)$, as well as their categories (total and non-continued).

Table 1. Members of the array $\boldsymbol{B}_\beta\left(\frac{7}{12},\frac{24}{41}\right)$ and their categories

| Denominators | Fractions | Multiples | Total Categorization | Non-continued Categorization | Label |
|---|---|---|---|---|---|
| 12 | 7/12 | | 3 | c | $\langle 3,\mathrm{c}\rangle$ |
| 13 | 8/13 | | 20 | 15 | $\langle 20,15\rangle$ |
| 14 | 8/14 | 4/7 | 15\|2 | sc | $\langle 15\|2,\mathrm{sc}\rangle$ |
| 15 | 9/15 | 3/5 | 16\|3 | c | $\langle 16\|3,\mathrm{c}\rangle$ |
| 16 | 9/16 | | 18 | 13 | $\langle 18,13\rangle$ |
| 17 | 10/17 | | 4 | 1 | $\langle 4,1\rangle$ |
| 18 | 11/18 | | 19 | 14 | $\langle 19,14\rangle$ |
| 19 | 11/19 | | 8 | 5 | $\langle 8,5\rangle$ |
| 20 | 12/20 | 3/5 | 16\|4 | c | $\langle 16\|4,\mathrm{c}\rangle$ |
| 21 | 12/21 | 4/7 | 15\|3 | sc | $\langle 15\|3,\mathrm{sc}\rangle$ |
| 22 | 13/22 | | 7 | 4 | $\langle 7,4\rangle$ |
| 23 | 13/23 | | 17 | 12 | $\langle 17,12\rangle$ |
| 24 | 14/24 | 7/12 | 3\|2 | c | $\langle 3\|2,\mathrm{c}\rangle$ |
| 25 | 15/25 | 3/5 | 16\|5 | c | $\langle 16\|5,\mathrm{c}\rangle$ |
| 26 | 15/26 | | 10 | 7 | $\langle 10,7\rangle$ |
| 27 | 16/27 | | 9 | 6 | $\langle 9,6\rangle$ |
| 28 | 16/28 | 4/7 | 15\|4 | sc | $\langle 15\|4,\mathrm{sc}\rangle$ |
| 29 | 17/29 | | 2 | sc | $\langle 2,\mathrm{sc}\rangle$ |

| Denominators | Fractions | Multiples | Total Categorization | Non-continued Categorization | Label |
|---|---|---|---|---|---|
| 30 | 18/30 | 3/5 | 16\|6 | c | ⟨16\|6, c⟩ |
| 31 | 18/31 | | 5 | 2 | ⟨5,2⟩ |
| 32 | 19/32 | | 11 | 8 | ⟨11,8⟩ |
| 33 | 19/33 | | 12 | 9 | ⟨12,9⟩ |
| 34 | 20/34 | 10/17 | 4\|2 | 1\|2 | ⟨4\|2,1⟩ |
| 35 | 20/35 | 4/7 | 15\|5 | sc | ⟨15\|5, sc⟩ |
| 36 | 21/36 | 7/12 | 3\|3 | c | ⟨3\|3, c⟩ |
| 37 | 22/37 | | 13 | 10 | ⟨13,10⟩ |
| 38 | 22/38 | 11/19 | 8\|2 | 5\|2 | ⟨8\|2,5⟩ |
| 39 | 23/39 | | 6 | 3 | ⟨6,3⟩ |
| 40 | 23/40 | | 14 | 11 | ⟨14,11⟩ |
| 41 | 24/41 | | 1 | c | ⟨1, c⟩ |

If we observe total categorization, it can be noticed that sorted array $\boldsymbol{B}_\beta$, with the least absolute error from the constant $\beta$ looks like this

$$\boldsymbol{\mathcal{B}}_\beta^\tau\left(\frac{7}{12},\frac{24}{41}\right)=\left(\frac{24}{41}\langle 1\rangle,\frac{17}{29}\langle 2\rangle,\frac{7}{12}\langle 3\rangle,\frac{14}{24}\langle 3|2\rangle,\frac{21}{36}\langle 3|3\rangle\,\frac{10}{17}\langle 4\rangle,\frac{20}{34}\langle 4|2\rangle,\right.$$
$$\frac{20}{34}\langle 4|2\rangle,\frac{18}{31}\langle 5\rangle,\frac{23}{39}\langle 6\rangle,\frac{13}{22}\langle 7\rangle,\frac{11}{19}\langle 8\rangle,\frac{22}{38}\langle 8|2\rangle,\frac{16}{27}\langle 9\rangle, \tag{39}$$
$$\frac{15}{26}\langle 10\rangle,\frac{19}{32}\langle 11\rangle,\frac{19}{33}\langle 12\rangle,\frac{22}{37}\langle 13\rangle,\frac{12}{20}\langle 16|4\rangle,\frac{15}{25}\langle 16|5\rangle,$$
$$\left.\frac{18}{30}\langle 16|6\rangle,\frac{13}{23}\langle 17\rangle,\frac{9}{16}\langle 18\rangle,\frac{11}{18}\langle 19\rangle,\frac{8}{13}\langle 20\rangle\right).$$

Let us notice that in previously shown record $\boldsymbol{\mathcal{B}}_\beta^\tau$ value written in angle brackets represents total category. If fraction is multiple after vertical line follows which number of multiple it is.

If we observe non-continued categorization, it can be noticed that sorted array $\boldsymbol{B}_\beta^\upsilon$ (Eq. 38), with the least absolute error from the constant $\beta$, looks like this

$$\boldsymbol{\mathcal{B}}_\beta^\upsilon\left(\frac{7}{12},\frac{24}{41}\right)=\left(\frac{10}{17}\langle 1\rangle,\frac{20}{34}\langle 1|2\rangle,\frac{18}{31}\langle 2\rangle,\frac{23}{39}\langle 3\rangle,\frac{13}{22}\langle 4\rangle\,,\frac{11}{19}\langle 5\rangle,\right.$$
$$\frac{22}{38}\langle 5|2\rangle,\frac{16}{27}\langle 6\rangle,\frac{15}{26}\langle 7\rangle,\frac{19}{32}\langle 8\rangle,\frac{19}{33}\langle 9\rangle,\frac{22}{37}\langle 10\rangle, \tag{40}$$

$$\frac{23}{40}\langle 11\rangle, \frac{13}{23}\langle 12\rangle, \frac{9}{16}\langle 13\rangle, \frac{11}{18}\langle 14\rangle, \frac{8}{13}\langle 15\rangle\Big)$$

Let us notice that in previously shown record $\boldsymbol{\mathcal{B}}_{\beta}^{\boldsymbol{v}}$ value written in angle brackets represents non-continued category. If fraction is multiple, after vertical line follows which number of multiple it is.

Let us emphasize that in the array $\boldsymbol{\mathcal{B}}_{\beta}^{\boldsymbol{v}}$ the first best non-continued fraction is $\frac{10}{17}$, is followed by his first multiple $\frac{20}{34}$, then the second best non-continued fraction is $\frac{18}{31}$, and so on. The best ranked non-continued fractions will be of importance in further analysis of the application of non-continued relational approximations.

Now, we can sort members of $\boldsymbol{B}_{\beta}$ by both categories:

$$\begin{aligned}\boldsymbol{\mathcal{B}}_{\beta}\left(\frac{7}{12}, \frac{24}{41}\right) = \Big(&\frac{24}{41}\langle 1, \mathrm{c}\rangle, \frac{17}{29}\langle 2, \mathrm{sc}\rangle, \frac{7}{12}\langle 3, \mathrm{c}\rangle, \frac{14}{24}\langle 3|2, \mathrm{c}\rangle, \frac{21}{36}\langle 3|3, \mathrm{c}\rangle,\\ &\frac{10}{17}\langle 4,1\rangle, \frac{20}{34}\langle 4|2,1\rangle, \frac{18}{31}\langle 5,2\rangle, \frac{23}{39}\langle 6,3\rangle, \frac{13}{22}\langle 7,4\rangle,\\ &\frac{11}{19}\langle 8,5\rangle, \frac{22}{38}\langle 8|2,5\rangle, \frac{16}{27}\langle 9,6\rangle, \frac{15}{26}\langle 10,7\rangle, \frac{19}{32}\langle 11,8\rangle,\\ &\frac{19}{33}\langle 12,9\rangle, \frac{22}{37}\langle 13,10\rangle, \frac{23}{40}\langle 14,11\rangle, \frac{8}{14}\langle 15|2, \mathrm{sc}\rangle,\\ &\frac{12}{21}\langle 15|3, \mathrm{sc}\rangle, \frac{16}{28}\langle 15|4, \mathrm{sc}\rangle, \frac{20}{35}\langle 15|5, \mathrm{sc}\rangle, \frac{9}{15}\langle 16|3, \mathrm{c}\rangle,\\ &\frac{12}{20}\langle 16|4, \mathrm{c}\rangle, \frac{15}{25}\langle 16|5, \mathrm{c}\rangle, \frac{18}{30}\langle 16|6, \mathrm{c}\rangle,\\ &\frac{13}{23}\langle 17,12\rangle, \frac{9}{16}\langle 18,13\rangle, \frac{11}{18}\langle 19,14\rangle, \frac{8}{13}\langle 20,15\rangle\Big).\end{aligned} \tag{41}$$

In the array $\boldsymbol{\mathcal{B}}_{\beta}$ in angle brackets we have three types of information, which are shown in last column of Table 1:

(1) Total category from the array $\boldsymbol{\mathcal{B}}_{\beta}^{\boldsymbol{\tau}}$,

(2) If the observe fraction is multiple, after vertical line it is specified which number of multiple it is, and

(3) After comma, continued fractions are specified by mark (c), semi-continued fractions by mark (sc), and non-continued fractions are specified by the value of position of the fraction based on non-continued categorization from the array $\boldsymbol{\mathcal{B}}_{\beta}^{\boldsymbol{v}}$.

We can notice that some multiples of continued and semi-continued fractions, such as $\frac{8}{14}$ (sc), $\frac{9}{15}$ (c), ... are worse approximate of constant β then the most of non-continued fractions from the array $\boldsymbol{\mathcal{B}}_{\beta}$, Table 1.

**Software Realization**

In order to determine the arrays $\boldsymbol{\mathcal{B}}^{\boldsymbol{\tau}}$, $\boldsymbol{\mathcal{B}}^{\boldsymbol{v}}$ and $\boldsymbol{\mathcal{B}}$ in general, we made a program in the programming language MatLab. Starting from the fact that, for the observed value α, we know continued and semi-continued fractions and that we determined them according to one of the methods from papers (Cortzen A. 2011; Malešević B. 1988; Malešević B, Milinković L. 2014.) we have started to determine the rest of fractions $\frac{k_i}{m_i}$. Let us denote that algorithms that are listed below are written in pseudo code.

**Algorithm 1.** Pseudo code for determining fractions $\frac{k_i}{m_i}$ and absolute error $\Delta_i$

```
function initialization
input α, N                // for α we are looking for first N best relational approximations
output A[]

for i=1 to N
    A[i].m = i
    temp = α· A[i].m
    A[i].k = round(temp)                    // Eq. 12
    A[i].Δ = abs(α –A[i].k/A[i].m)
end for
```

**Algorithm 2.** Pseudo code for categorization

```
function categorization
input B[]
```

```
output B_out[]

B_out = sort(B, growing by B.Δ)

cat = 0                              // tau means total category
for i=1 to length(B.k)
    B_out[i].multipl= gcd(B_out[i].k, B_out[i].m)
    if (B_out[i].multipl == 1) or (i==1) then
        cat = cat + 1
    else if B_out[i].m/B_out[i].multipl ≠ B_out[i-1].m/B_out[i-1].multipl then
            cat = cat + 1
    end if
    B_out[i].cat = cat
end for
```

**Algorithm 3.** Pseudo code for determining the arrays $\boldsymbol{\mathcal{B}^{\tau}}$, $\boldsymbol{\mathcal{B}^{\upsilon}}$ and $\boldsymbol{\mathcal{B}}$

```
function classification
input B[], B^c[], B^sc[]                 // Eq. 31
output 𝓑^τ[], 𝓑^υ[]

𝓑^τ = categorization(B)

C = B – ((B ∩ B^c) + (B ∩ B^sc))
𝓑^υ = categorization(C)

for i=1 to length(𝓑^τ)
    𝓑[i].k = 𝓑^τ[i].k
    𝓑[i].m = 𝓑^τ[i].m
    𝓑[i].tau = 𝓑^τ[i].cat
    𝓑[i].multipl = 𝓑^τ[i].multipl
    temp = seek(𝓑[i], 𝓑^υ)               // Function seeks position of element 𝓑[i] in array
                                         // of non-continued fractions 𝓑^υ.
                                         // If element is not found returns empty.
    if (temp ≠ empty) then
        𝓑[i].ni = 𝓑^υ[temp].cat
    else
        temp = seek(𝓑[i], B^c)           // Function seeks position of element 𝓑[i] in array
                                         // of continued fractions B^c.
        if (temp ≠ empty) then
            𝓑[i].ni = c
        else
            𝓑[i].ni = sc
        end if
    end if
end for
```

## Problem of the Microtonal Scale

During research of the microtonal scales we have noticed that rational approximation $\frac{10}{17}$ was used as a good approximation of the music constant, (according to Benson D. 2006; Carey N, Clampitt D. 1989; Douthett J, Krantz R. 2007; Khrushchev S. 2008; Noll T. 2006; Schechter M. 1980). This value determines microtonal scale 17EDO. Based on our analysis in this paper, it is mathematically clear that $\frac{10}{17}$ is not the best rational approximation of the first type, and also, that it is not continued or semi-continued fraction. Let us notice that in fraction array from $\frac{7}{12}$ to $\frac{24}{41}$ the best non-continued rational approximation of the constant $\beta$ is fraction $\frac{10}{17}$, (Eq. 40).

As a good approximation of the musical constant β can be used fraction $\frac{18}{31}$, (Benson D. 2006). This is also a non-continued fraction and this value determines microtonal scale 31EDO. Let us point out that this microtonal scale represents Dutch national scale choice, of which is written in websites (Fokker A. D; Terpstra S; 31 Tone Equal Temperament) as well as in works (Dillon G, Musenich R. 2009; Keislar D. 1991; Khramov M. 2009; Rapoport P. 1987). In analysed fraction array from $\frac{7}{12}$ to $\frac{24}{41}$ before fraction $\frac{18}{31}$ there are only two non-continued fractions that are better approximations of musical constant $\beta$ (fraction $\frac{10}{17}$ and its multiple $\frac{20}{34}$). Thus, fraction $\frac{18}{31}$ is the second best non-continued approximation of the musical constant in considered range, (Eq. 40).

## Conclusion

This paper introduces and considers non-continued approximations. These allow categorization of microtonal scales by the quality of approximation of musical constant $\beta$ that defines perfect fifth. In particular, the quality of the two frequently used microtonal scales 17EDO and 31EDO were estimated. The fractions which are corresponding to these microtonal scales $\left(\frac{10}{17}\text{ and }\frac{18}{31}\right)$ are also the best two non-continued approximations in fraction array from $\frac{7}{12}$ to $\frac{24}{41}$ for the constant $\beta$.

**Acknowledgement.**

This research was supported by the Serbian Ministry of Education and Science (projects no. ON 174032, III 44006 and TR 36026).